\documentstyle[12pt,amssymb]{article}
\begin{document}
\newtheorem{Theo}{Theorem}
\newtheorem{Ex}{Example}
\newtheorem{Def}{Definition}
\newtheorem{Lem}{Lemma}
\newtheorem{Cor}{Corollary}
\newenvironment{pf}{{\noindent\bf Proof.\ }}{\ $\Box$\medskip}
\newtheorem{Prop}{Proposition}
\renewcommand{\theequation}{\thesection.\arabic{equation}}


\mathchardef\za="710B  
\mathchardef\zb="710C  
\mathchardef\zg="710D  
\mathchardef\zd="710E  
\mathchardef\zve="710F 
\mathchardef\zz="7110  
\mathchardef\zh="7111  
\mathchardef\zvy="7112 
\mathchardef\zi="7113  
\mathchardef\zk="7114  
\mathchardef\zl="7115  
\mathchardef\zm="7116  
\mathchardef\zn="7117  
\mathchardef\zx="7118  
\mathchardef\zp="7119  
\mathchardef\zr="711A  
\mathchardef\zs="711B  
\mathchardef\zt="711C  
\mathchardef\zu="711D  
\mathchardef\zvf="711E 
\mathchardef\zq="711F  
\mathchardef\zc="7120  
\mathchardef\zw="7121  
\mathchardef\ze="7122  
\mathchardef\zy="7123  
\mathchardef\zf="7124  
\mathchardef\zvr="7125 
\mathchardef\zvs="7126 
\mathchardef\zf="7127  
\mathchardef\zG="7000  
\mathchardef\zD="7001  
\mathchardef\zY="7002  
\mathchardef\zL="7003  
\mathchardef\zX="7004  
\mathchardef\zP="7005  
\mathchardef\zS="7006  
\mathchardef\zU="7007  
\mathchardef\zF="7008  
\mathchardef\zW="700A  

\newcommand{\be}{\begin{equation}}
\newcommand{\ee}{\end{equation}}
\newcommand{\lra}{\longrightarrow}
\newcommand{\ra}{\rightarrow}
\newcommand{\bea}{\begin{eqnarray}}
\newcommand{\eea}{\end{eqnarray}}
\newcommand{\beas}{\begin{eqnarray*}}
\newcommand{\eeas}{\end{eqnarray*}}
\newcommand{\Z}{{\Bbb Z}}
\newcommand{\R}{{\Bbb R}}
\newcommand{\C}{{\Bbb C}}
\newcommand{\SL}{SL(2,\R)}
\newcommand{\Sl}{sl(2,\C)}
\newcommand{\SU}{SU(2)}
\newcommand{\su}{su(2)}
\newcommand{\G}{{\cal G}}
\newcommand{\g}{{\frak g}}
\newcommand{\h}{{\frak h}}
\newcommand{\D}{{\rm d}}
\newcommand{\de}{\,{\stackrel{\rm def}{=}}\,}
\newcommand{\Ad}{{\rm Ad}}
\newcommand{\we}{\wedge}
\newcommand{\We}{\bigwedge}
\newcommand{\tL}{{\tilde\zL}}
\newcommand{\tl}{{\tilde\zl}}
\newcommand{\ta}{{\tilde\za}}
\newcommand{\tda}{{\widetilde{\D\za}}}
\newcommand{\nn}{\nonumber}
\newcommand{\ot}{\otimes}
\newcommand{\s}{{\textstyle *}}
\newcommand{\Li}{{\cal L}}
\newcommand{\const}{{\rm const}}
\newcommand{\pa}{\partial}
\newcommand{\ti}{\times}
\begin{center}
{\Large\bf Remarks on Nambu-Poisson,\\ and Nambu-Jacobi brackets}
\vskip 1cm
J. Grabowski\footnote{Institute of Mathematics, Warsaw University,
ul. Banacha 2, 02-097 Warszawa, Poland;
{\it e-mail:} jagrab@mimuw.edu.pl .\\
This work has been supported by KBN, grant No. 2 P03A 042
10.} \\ and \\
G. Marmo\footnote{Dipartimento di Scienze Fisiche, Universit\`a di Napoli,
Mostra d'Oltremare, Pad. 20, 80125 Napoli, Italy;
{\it e-mail:} gimarmo@na.infn.it .\\
This work has been partially supported by PRIN-97 "SINTESI".}
\end{center}
\date{\ }
\centerline{\bf Abstract}
We show that Nambu-Poisson and Nambu-Jacobi brackets can be defined
inductively:  an  $n$-bracket,  $n>2$,  is   Nambu-Poisson   (resp.
Nambu-Jacobi) if  and  only  if  fixing  an  argument  we  get  an
$(n-1)$-Nambu-Poisson  (resp.   Nambu-Jacobi)   bracket.   As    a
by-product we get relatively simple proofs   of Darboux-type
theorems for these structures.
\bigskip\noindent
\setcounter{equation}{0}
\section{Introduction}
The concept of a Nambu-Poisson structure was introduced by Takhtajan
\cite{Ta} in order to find an axiomatic formalism for
 the $n$-bracket operation
\be\label{01}
\{ f_1,\dots,f_n\}={\rm det}\left(\frac{\pa f_i}{\pa x_j}\right),
\ee
proposed  by  Nambu  \cite{Nam}  to  generalize  the   Hamiltonian
mechanics (cf. also \cite{BF, Cha, FDS}).
Takhtajan \cite{Ta} has observed that the Nambu canonical  bracket
(\ref{01}) is $n$-linear skew-symmetric and satisfies the
{\em generalized Jacobi identity}:
\bea\label{J}
&\{ f_1,\dots,f_{n-1},\{ g_1,\dots,g_n\}\}=
\{\{ f_1,\dots,f_{n-1},g_1\},g_2,\dots,g_n\}+\\
&\{ g_1,\{ f_1,\dots,f_{n-1},g_2\},g_3,\dots,g_n\}+\dots+
\{ g_1,\dots,g_{n-1},\{ f_1,\dots,f_{n-1},g_n\}\}.\nn
\eea
Such an axiom was also considered by other authors about the  same
time (see \cite{SV}).
These, however, are exactly the axioms of an $n$-Lie algebra introduced
by Filippov \cite{Fi} in 1985, who gave also the example of the canonical
Nambu bracket (\ref{01}) in this context. The additional assumption made
by Takhtajan was that the bracket, acting on smooth functions, has to
satisfy the Leibniz rule
\be\label{L}
\{ fg,f_2,\dots,f_{n-1}\}=f\{ g,f_2,\dots,f_{n-1}\}+
\{ f,f_2,\dots,f_{n-1}\} g,
\ee
what generalizes the notion of a Poisson bracket and means that the
bracket is in fact defined by an $n$-vector field $\zL$ in the obvious
way:
\be
\{ f_1,\dots,f_n\}=\zL_{f_1,\dots,f_n},
\ee
where we denote $\zL_{f_1,\dots,f_k}$ to be the contraction
$i_{\D f_k}\cdots i_{\D f_1}\zL$.
The generalized Jacobi identity (\ref{J}) means then that the
{\em hamiltonian vector fields} $\zL_{f_1,\dots,f_{n-1}}$
(of $(n-1)$-tuples of functions
this time)  preserve the tensor $\zL$,
i.e. the corresponding Lie derivative
(which we write as the Schouten bracket) vanishes:
\be\label{LD}
[\zL_{f_1,\dots,f_{n-1}},\zL]=0.
\ee
This implies also that the characteristic distribution $D_\zL$ of
the $n$-vector field $\zL$, i.e. the distribution generated by
all the hamiltonian vector fields, is involutive.
Indeed, from (\ref{J}) we easily derive
\be\label{ham}
[\zL_{f_1,\dots,f_{n-1}},\zL_{g_1,\dots,g_{n-1}}]=
\sum_i\zL_{g_1,\dots,\{ f_1,\dots,f_{n-1},g_i\},\dots,g_{n-1}}.
\ee
We have even more: the characteristic distribution defines a (possible
singular) foliation of the manifold $M$ in the sense of Stefan-Sussmann.
It is due to the fact that the module of 1-forms on $M$ is finitely
generated over the ring $C^\infty (M)$ of smooth functions on $M$,
so we can generate the distribution by a finite number of hamiltonian
vector fields which are closed over $C^\infty (M)$ under the Lie
bracket and the Stefan-Sussmann condition of integrability of the
distribution is satisfied. The proof goes exactly parallel to that
in the classical Poisson case.
All this look quite similar to the case of classical Poisson structures.
Now, the point is that in the case of Nambu-Poisson structures of order
$n>2$ the leaves of the characteristic foliation have to be either 0 or
$n$-dimensional; a different behavior comparing with the classical
Poisson case. We shall consider this point later.
\par
Since Nambu-Poisson brackets are just the Filippov brackets which are
given by multi-derivations of the associative algebra $C^\infty (M)$
(the Leibniz rule), the next obvious generalization is to follow the
idea of Kirillov \cite{Ki} and to ask what are Filippov $n$-brackets on
the ring of functions, which are given by local operators. They could
be called {\em Nambu-Jacobi} brackets since the Jacobi brackets are
what we get in the binary case (cf. \cite{Gr} for a purely algebraic
approach).
It is easy to see from (\ref{J}) that every contraction of the bracket
with an element $f$ leads from an $n$-Lie algebra in the sense of Filippov
(let us call them simply {\em $n$-Filippov algebras}) to an $(n-1)$-Filippov
algebra: the $(n-1)$-ary bracket $[\ ,\dots,\ ]_f$ defined by
$[f_1,\dots,f_{n-1}]_f=[f,f_1,\dots,f_{n-1}]$ satisfies the $(n-1)$-Jacobi
identity if the bracket $[\ ,\dots \ ]$ satisfies the $n$-Jacobi identity.
Hence, contractions lead from $n$-Nambu-Jacobi brackets to $(n-1)$-
Nambu-Jacobi brackets (we will show later that it can be inverted).
The binary Jacobi brackets are given by first-order differential operators
(\cite{Ki,Gr}), so, due to the fact that they are totally skew-symmetric,
$n$-Nambu-Jacobi bracket are also given by first-order differential operators.
Similarly as in the binary case, they can be written with the help of two
multivector fields: $n$-vector field $\zD$ and $(n-1)$-vector field $\zG$,
in the form
\be\label{02}
\{ f_1,\dots,f_n\}=(\zD+s(\zG))(f_1,\dots,f_n),
\ee
where $\zD(f_1,\dots,f_n)=\zD_{f_1,\dots,f_n}$ is just the bracket
induced by $\zD$ and
\be
s(\zG)(f_1,\dots,f_n)=\sum_i(-1)^{i+1}f_i\zG_{f_1,\dots,f_{i-1},f_{i+1},
\dots,f_n}
\ee
(cf. \cite{MVV}).
The $n$-Jacobi identity puts additional restrictions for the pair
$(\zD,\zG)$ to get a Nambu-Jacobi structure. For the classical case,
$n=2$, they read
\bea
&[\zG,\zD]=0, \label{x}\\
&[\zD,\zD]=-2\zG\we\zD ,\label{xx}
\eea
where  the  brackets  are   the  Schouten  brackets.
We use the Schouten bracket
\bea [X_1\wedge\dots \wedge X_k, Y_1\wedge \dots
\wedge Y_k] = \sum_{i,j}(-1)^{i+j} [X_i,Y_j] \wedge X_1\wedge\\
\dots \wedge \widehat{X}_i \wedge \dots \wedge X_k\wedge
Y_1\wedge \dots \wedge \widehat{Y}_j\wedge \dots \wedge Y_l.\nn
\eea
Let us note that sometimes one can meet a  version  of  (\ref{xx})
differing by sign (cf. \cite{Li,DLM}) when the bracket  differs  by
sign from the bracket we use.
\par
In this paper we prove the  inductive
property of Nambu-Poisson and Nambu-Jacobi bracket: for $n>2$, an
$n$-linear skew-symmetric
bracket of smooth functions  on  a  manifold  is  a  Nambu-Poisson
(resp. Nambu-Jacobi) bracket if and only if fixing one argument we
get  a Nambu-Poisson (resp. Nambu-Jacobi) bracket of order $(n-1)$.
As a by-product we get versions of Darboux-type  theorems  for
these bracket (cf. \cite{AG,Ga,Na1,Pa,MVV}) with relatively  short
proofs.
\par
There are other concepts  of  $n$-ary  Lie,  Poisson,  and  Jacobi
brackets using
a generalized Jacobi identity of different type than (\ref{J}), a
skew-symmetrization of it. We will not discuss them here, so let
us  only  mention  the  papers   \cite{APP,AIP,ILMD,ILMP,MV}   and
references there. Recently a unification of different concepts was
proposed in \cite{VV}.

\section{Recursive   characterization   of    Nambu-Poisson    and
Nambu-Jacobi algebras}
We start with the following easy observation.
\begin{Prop} Let $X_1,\dots,X_n$ be vector fields on a manifold $M$.
Then $\zL=X_1\we\dots\we X_n$ is a Nambu-Poisson tensor if and only
if the distribution $D$ generated by these vector fields is involutive
at regular points of $\zL$.
\end{Prop}
\begin{pf} The proposition has a local character, so all we shall do
will be local near regular points and we may assume that the vector
fields are linearly independent. Under this condition $D$ coincides with the
characteristic distribution $D_\zL$ of the tensor field $\zL$, i.e.
\be
D={\rm span}\{\zL_{f_1,\dots,f_{n-1}}: f_i\in C^\infty (M), i=1,
\dots,n-1 \}.
\ee
If $\zL$ is a Nambu-Poisson tensor then $D_\zL=D$ is known to be involutive.
On the other hand, if $D$ is involutive, then it generates an $n$-dimensional
foliation which, in appropriate coordinates, is generated by the coordinate
vector fields $\pa_1,\dots,\pa_n$ ( a version of the Frobenius theorem), so
that $\zL=f\pa_1\we\dots\we\pa_n$ with some function $f$.
This is a standard example of a Nambu-Poisson tensor (cf. \cite{MVV},
Corollary
3.2). As a matter of fact, the function $f$ could be chosen to be a constant.
\end{pf}

\medskip\noindent
{\bf Remark.} At singular points of $\zL$ the situation may be different.
For instance, $\zL=\pa_1\we(x_1\pa_2)$ is a Poisson tensor, but $D$ is
not involutive at points, where $x_1=0$: $[\pa_1,x_1\pa_2]=\pa_2\notin D$.

\medskip\noindent
We shall make use of the following variant of the lemma `on three planes'
(cf. \cite{MVV} or \cite{DZ}).
\begin{Lem}
Let $\{\zL_i: i\in I\}$  be  a  family  of  decomposable  non-zero
$n$-vectors  of  a  vector  space  $V$   such   that   every   sum
$\zL_{i_1}+\zL_{i_2}$ is again decomposable.
Then,
\item{(a)} the linear span $D$ of the linear subspaces  $D_{\zL_i}$
they generate is at most $(n+1)$-dimensional
\item{}or
\item{(b)} the  intersection $\cap_iD_{\zL_i}$ is at least
$(n-1)$-dimensional.
\end{Lem}
\begin{pf}
It  is  easy  to  see  that  the  sum   $\zL_{i_1}+\zL_{i_2}$   is
decomposable, where the summands are non-zero, if and only if the
intersection of $n$-dimensional subspaces
$D_{\zL_{i_1}}\cap D_{\zL_{i_2}}$ is at least $(n-1)$-dimensional.
Then we can use a corrected version of `lemma on three planes'
(\cite{MVV}, Lemma 4.4.) as in \cite{DZ},  which  states  that  in
this case we have (a) or (b), with rather obvious proof.
\end{pf}

\medskip\noindent
\begin{Lem}
Let $\zL$ be an n-vector field on a manifold $M$, $n>2$, such that
all the contractions $\zL_f$, with $f\in C^\infty (M)$, are Nambu-Poisson
tensors. Then $\zL$ is decomposable at its regular points.
\end{Lem}
\begin{pf} The Nambu-Jacobi identity for $\zL_f$ reads
\be
[\zL_{f,f_1,\dots,f_{n-2}},\zL_f]=0,
\ee
where $[\ ,\ ]$ stands for the Schouten bracket. It easily implies that
the operation $A$ acting on functions by
\be
A(f_1,\dots,f_{n-1},g_1,\dots,g_n)=[\zL_{f_1,\dots,f_{n-1}},\zL]_
{g_1,\dots,g_n}
\ee
is totally skew-symmetric (being skew-symmetric with respect to $f_i$'s,
$g_i$'s and vanishing for $f_1=g_1$) and it is represented by a $(2n-1)$-
vector field, since it acts by derivations on $g_i$'s.
This implies that
\be
A(f_1^2,f_2,\dots,f_{n-1},g_1,\dots,g_n)=2f_1A(f_1,f_2,\dots,f_{n-1},
g_1,\dots,g_n).
\ee
On the other hand, from properties of the Schouten bracket we get
\bea\nn
&A(f_1^2,f_2,\dots,f_{n-1},g_1,\dots,g_n)=
[2f_1\zL_{f_1,f_2,\dots,f_{n-1}},\zL]_{g_1,\dots,g_n}=\\
&2f_1A(f_1,f_2,\dots,f_{n-1},g_1,\dots,g_n)\pm 2
(\zL_{f_1,f_2,\dots,f_{n-1}}\we\zL_{f_1})_{g_1,\dots,g_n},
\eea
so that
\be
\zL_{f_1,f_2,\dots,f_{n-1}}\we\zL_{f_1}=0.
\ee
The last identity, for $n>2$, implies that $\zL$ is decomposable at regular
points, as shown in \cite{MVV}, Proposition 4.1, or \cite{Ga}.
\end{pf}
\medskip\noindent
\begin{Theo} An $n$-vector field $\zL$, $n>2$, on a manifold $M$ is a
Nambu-Poisson tensor if and only if its contractions $\zL_f$ are Nambu-
Poisson tensors for all $f\in C^\infty (M)$.  In  this  case,  the
tensor $\zL$ can be written in appropriate coordinates around  its
regular points in the form $\pa_1\we\dots\we\pa_n$.
\end{Theo}
\begin{pf}
The implication $\Rightarrow$ is well known and trivial.  To  show
the converse, we can first reduce to regular points and then to
make use of Lemma 2 to get that $\zL$ is decomposable, say
$\zL=X_1\we\dots\we X_n$.
Now, we have to show that the distribution $D$ generated by the linearly
independent vector fields $X_1,\dots,X_n$ is involutive.
Since $X_1$ is (locally) non-vanishing,
we can find a function $f$ such that (again locally) $X_1(f)\equiv 1$.
Putting now $X'_i=X_i-X_i(f)X_1$ for $i>1$, we get
$\zL=X_1\we X'_2\we\dots\we X'_n$ with $X'_i(f)=0$.
Thus $\zL_f=X'_2\we\dots\we X'_n$ is a Nambu-Poisson tensor by the
inductive assumption, so that, in certain local coordinates,
$\zL_f=\pa_1\we\dots\we\pa_{n-1}$ as in the proof of Proposition1, and
\be
\zL=X_1\we\zL_f=X_1\we\pa_1\we\dots\we\pa_{n-1}.
\ee
In fact, we can additionally assume that $X_1(x_i)=0$ for $i=1,\dots,n-1$,
replacing $X_1$ by $X_1-\sum_1^{n-1}X_1(x_i)\pa_i$.
Assuming that the characteristic distribution $D$ generated by $\zL$ is not
involutive, we would get that $[\pa_i,X_1]\notin D$ for some $i$, say
$i=n-1$. But then $\zL_{x_1}=-X_1\we\pa_2\we\dots\we\pa_{n-1}$ generates
a distribution which is not involutive, contrary to the inductive assumption.
\end{pf}
\medskip\noindent
\begin{Cor} An $n$-linear
skew-symmetric bracket $\{\ ,\dots,\ \}$ on functions
on a manifold $M$, $n>2$, is a Nambu-Poisson bracket if and only if its
contraction with any $(n-2)$ functions gives a Poisson bracket.
\end{Cor}
\begin{pf} Since the contractions give Poisson bracket, our bracket
operation is given by an $n$-linear first-order differential operator
vanishing on constants, so by an $n$-vector field $\zL$. The rest
follows by applying Theorem 1 recursively.
\end{pf}

\medskip\noindent
We have a similar theorem for Nambu-Jacobi brackets.
They have to be of the form $\zD+s(\zG)$ (cf. \cite{MVV}) and we get
particular cases: just the Nambu-Poisson structure $\zD$ (locally $\zG=0$)
and the bracket given by $s(\zG)$ (locally, $\zD=0$).
Other examples we get putting in local coordinates $\zD=\pa_1\we\dots\we\pa_n$
and $\zG=\pa_1\we\dots\we\pa_{n-1}$. The characterization theorem (\cite{MVV})
shows that this is rather general picture. The following inductive theorem
will suggest an alternative proof of this characterization theorem.

\begin{Theo} An $n$-linear
skew-symmetric bracket $\{\ ,\dots,\ \}$ on functions
on a manifold $M$, $n>2$, is a Nambu-Jacobi bracket if and only if its
contraction $\{ f,\ ,\dots,\ \}$
with any function $f$ is a Nambu-Jacobi bracket.
Moreover,
\be
\{ f_1,\dots,f_n\}=\zD_{f_1,\dots,f_n}+s(\zG)(f_1,\dots,f_n)
\ee
for some multivector fields $\zD$ and $\zG$ which, in some local coordinates
around points where they do not vanish, can be written in the form
$\zD=\pa_1\we\dots\we\pa_n$ and $\zG=\pa_1\we\dots\we\pa_{n-1}$.
If, locally, one of the  tensors  $\zD,\zG$  vanishes,  the
other is a Nambu-Poisson tensor  ($\zG$  is  an  ordinary  Poisson
tensor of rank 2 if $n=3$) and can be written as above around  its
regular points.
\end{Theo}
\begin{pf} Only the implication $\Leftarrow$ is non-trivial.
Since the contractions are given by first-order differential operators,
our bracket is given by an $n$-linear first-order differential operator,
i.e. (cf. \ref{02})
\be
\{ f_1,\dots,f_n\}=\zD_{f_1,\dots,f_n}+s(\zG)(f_1,\dots,f_n)
\ee
for an $n$-vector field $\zD$ and an $(n-1)$-vector field $\zG$.
It is easy to see that the analogous decomposition for the contraction
with the function $f$ yields
\be
\{ f,\cdot,\dots,\cdot\}=(\zD_f+f\zG)-s(\zG_f),
\ee
i.e. $(\zD_f+f\zG,-\zG_f)$ is a Nambu-Jacobi structure for all $f\in
C^\infty (M)$. In particular, for $f\equiv 1$, we get that $(\zG,0)$ is
a Nambu-Jacobi structure, which implies that $\zG$ is a Nambu-Poisson
tensor. If, locally, $\zD\equiv 0$, then our original bracket is of
the form $s(\zG)$, so it is a Nambu-Jacobi structure.
If, on the other hand, $\zG$ vanishes locally, then our original
bracket is just given by $\zD$. But the contractions $\zD_f$ give
Nambu-Jacobi and hence Nambu-Poisson ($\zG=0$) structures, so our
bracket is a Nambu-Poisson bracket in view of Theorem 1.
\par
Since our theorem is local and it suffices to check the Jacobi identities
on an open-dense subset, we may assume now that $\zD\ne 0$ and $\zG\ne 0$.
\par
1. The case $n=3$.
\par
Since for all $f\in C^\infty (M)$ the pairs $(\zD_f+f\zG,-\zG_f)$
constitute the usual Jacobi structures, according  to  (\ref{xx}),
we have the identity
\be
[\zD_f+f\zG,\zD_f+f\zG]=2\zG_f\we(\zD_f+f\zG).
\ee
Computing the Schouten brackets and using the fact that $\zG$ is an
ordinary Poisson structure ($[\zG,\zG]=0$ and $[\zG_f,\zG]=0$), we get
\be
[\zD_f,\zD_f]+2f[\zD_f,\zG]-2f\zG_f\we\zG=2\zG_f\we\zD_f+2f\zG_f\we\zG
\ee
(let us notice that $[\zG,f]=-\zG_f$) which gives
\be\label{1}
[\zD_f,\zD_f]+2f[\zD_f,\zG]-2\zG_f\we\zD_f-4f\zG_f\we\zG=0.
\ee
Putting $f:=f+1$ in (\ref{1}), we see that
\be\label{2}
[\zD_f,\zG]-2\zG_f\we\zD=0
\ee
and hence
\be\label{3}
[\zD_f,\zD_f]+2\zG_f\we\zD_f=0.
\ee
Further, replacing $f$ by $f^2$ in (\ref{2}), we get
\bea\label{4}
0=[\zD_{f^2},\zG]-2\zG_{f^2}\we\zD=[2f\zD_f,\zG]-4f\zG_f\we\zD=\\
2f[\zD_f,\zG]-2\zG_f\we\zD_f-4f\zG_f\we\zD=-2\zG_f\we\zD_f\nn
\eea
which, compared with (\ref{3}), gives
\be
[\zD_f,\zD_f]=0
\ee
for all $f\in C^\infty (M)$. The last  means  that  $\zD_f$  is  a
Poisson tensor for each $f\in C^\infty (M)$, so $\zD$ is itself a
Nambu-Poisson tensor (and hence decomposable) according to Theorem 1.
\par
Assume now that $\zD\ne 0$.
In view of (\ref{4}), $\zG_f$ divides the decomposable tensor
$\zD_f$, if only  $\zD_f\ne  0$,  and  hence  divides  also  $\zD$.
If $\zD_f=0$ at a given point, we can use the linearization
\be
\zG_f\we\zD_g+\zG_g\we\zD_f=0
\ee
of (\ref{4}) to get (at this point) $\zG_f\we\zD_g=0$, for  $g$
chosen  such  that
$\zD_g\ne 0$, and to conclude that $\zG_f\we\zD=0$ for all $f$.
\par
Now, similarly as in (\cite{MVV}, Theorem
5.1), we can find local coordinates such that $\zG=\pa_1\we\pa_2$ and
$\zD=\phi\pa_1\we\pa_2\we\pa_3$ for some function $\phi$.
Now, since (\ref{1}) implies that $[\zG_f,\zD_f]=0$, and putting
$f=x_i$, $i=1,2$, we get that $[\pa_i,\phi\pa_i\we\pa_3]=0$, $i=1,2$,
so that the vector fields $\pa_1,\pa_2,\phi\pa_3$ pairwise commute
and we can chose local coordinates so that $\phi\equiv 1$.
In particular, $\zG=\zD_{x_3}$.
\par
If, locally, $\zD=0$, then (\ref{1}) reduces  to  $f\zG_f\we\zG=0$
for all $f$ and hence $\zG$ is decomposable, i.e. $\zG$ is an
ordinary Poisson tensor of rank 2.
\par
2. The case $n>3$.
\par
Since $(\zD_f+f\zG,-\zG_f)$ is a Nambu-Jacobi structure of order
$(n-1)>2$, $\zD_f+f\zG$ is a decomposable (at its regular points)
Nambu-Poisson tensor.
Let $D_f$ be its characteristic distribution. The decomposable
tensor $\zD_{f+g}+(f+g)\zG$ is the sum of two decomposable tensors
$(\zD_f+f\zG)+(\zD_g+g\zG)$ at their regular points, so the dimension
of $D_f\cap D_g$ is at least $(n-2)$ if $D_f,D_g\ne\{ 0\}$.
It is easy to see that the intersection of non-trivial $D_f$  must
be zero.
\par
Indeed, if a vector field $X$ divides all $\zD_f+f\zG$, then $X$
divides $\zG$ (put $f\equiv 1$) and hence all  $\zD_f$,  so  that
$\zD_f=X\we  Y^f$, where
$Y^f$ is an $(n-2)$-vector field (all is local and for  "most"  of
functions $f$ we have $\zD_f\ne 0$).
Finding a function $g$ such that $X(g)=1$, we can assume, as in the
proof of Theorem 1, that $Y^f_g=0$, so that $\zD_{f,g}=Y^f\ne 0$.
Similarly,  $\zD_{g,g}=Y^g$,  but  $\zD_{g,g}=0$,  so  $\zD_g=X\we
Y^g=0$; a contradiction, since $\zD_{f,g}=-\zD_{g,f}\ne 0$.
Therefore, Lemma 1 implies that the linear span $D$ of all
the distributions $D_f$ has the dimension $\le n$. Moreover,
since $D_\zG=D_1$, also $D_\zG\subset D$. Now,
\be
(\zD_{f_1}+f_1\zG)_{f_2,\dots,f_{n-1}}=\zD_{f_1,f_2,\dots,f_{n-1}}+
f_1\zG_{f_2,\dots,f_{n-1}}\in D
\ee
implies that $\zD_{f_1,f_2,\dots,f_{n-1}}\in D$, i.e. the characteristic
distribution $D_\zD$ of $\zD$ is contained in $D$. But $\zD$ is an $n$-tensor,
so $dim(D_\zD)\ge n$ at regular points and hence $D_\zD=D$ and $\zD$ is
decomposable at its regular points.
\par
This implies that we can find local coordinates such that
$\zG=\pa_1\we\dots\we\pa_{n-1}$ and $\zD=\pa_1\we\dots\we\pa_{n-1}\we X$
for a vector field $X$ which we may assume to annihilate the coordinate
functions $x_1,\dots,x_{n-1}$.
If the characteristic distribution $D=D_\zD$ would be not involutive,
say $[\pa_{n-1},X]\notin D$ then the Nambu-Poisson tensor
\be
\zD_{x_1}+x_1\zG=\pa_2\we\dots\we\pa_{n-1}\we(X\pm x_1\pa_1)
\ee
would have a non-involutive characteristic distribution; a contradiction.
Therefore $\zD$ is a Nambu-Poisson tensor and we can write the vector
field $X$ in the form $X=\phi\pa_n$ with a function $\phi$.
Let us show now that the hamiltonian vector fields of $\zG$ preserve
$\zD$. Indeed, since $(\zD_f+f\zG, -\zG_f)$ is a Nambu-Jacobi structure,
at regular points $-\zG_f=(\zD_f+f\zG)_h$ for some function $h$ and
hence
\be
-[\zG_{f,f_1,\dots,f_{n-3}},\zD_f+f\zG]=
[(\zD_f+f\zG)_{h,f_1,\dots,f_{n-3}},\zD_f+f\zG]=0
\ee
($\zD_f+f\zG$ is a Nambu-Poisson tensor). But $[\zG_{f,f_1,\dots,
f_{n-3}},f\zG]=0$ ($\zG$ is a Nambu-Poisson tensor), so that
\be
[\zG_{f,f_1,\dots,f_{n-3}},\zD_f]=0.
\ee
Putting now $(f,f_1,\dots,f_{n-3})=(x_i,x_1,\dots,\check{x_i},\dots,
\check{x_j},\dots,x_{n-1}),$ where $\check{x_i}$ stands for omission,
we get, similarly as in the case $n=3$, that
$\pa_1,\dots,\pa_{n-1},\phi\pa_n$ are pairwise commuting vector fields,
so we can put $\phi\equiv 1$ in appropriate local coordinates.
In particular, $\zG=(-1)^n\zD_{x_n}$.
\end{pf}
\section{Conclusions}  In  last  few  years  interest  to  $n$-ary
generalizations of the  concept  of  Lie  algebra,  especially  to
Nambu-Poisson and Nambu-Jacobi brackets, has  been  growing  among
mathematicians and physicists.
\par
We   proved   that    Nambu-Poisson    and
Nambu-Jacobi bracket can be defined inductively: a bracket is such
an $n$-bracket if and only if, when contracted with any  function,
it gives such  an $(n-1)$-bracket.  This  reduces  the  general
case  to classical Poisson and Jacobi structures. In the proof we
have strongly used  the  fact  that  these  bracket  are  given  by
decomposable multivector fields defining first-order differential
operators  on the algebra of smooth functions on a manifold.
\par
The question, if an analog  of  this  fact  is  true  for   any
(finite-dimensional) Filippov algebra, remains open. Any Filippov
$n$-algebra  structure on a vector space $V$ is defined by a linear
multivector field
\be    \zL=\sum_{i_1,\dots,i_n}c^k_{i_1,\dots,i_n}x_k\pa_{x_{i_1}}
\we\cdots\we\pa_{x_{i_n}}
\ee
on $V^*$, where a basis $(x_i)$ of $V$ is regarded as a coordinate
system for $V^*$. This tensor  field  defines  an  $n$-bracket  on
linear functions  which  should  satisfy  the  generalized  Jacobi
identity. In  general,  however,  such  tensors  need  not  to  be
decomposable,  since  to  direct  products  of  Filippov  algebras
correspond `direct sums' of the  corresponding  tensors which  can
never be decomposable (if the summands are  non-zero).  Therefore,
we finish with the following problem.

\bigskip\noindent
{\bf  Problem.} {\it Let  $[\cdot,\dots,\cdot]$  be  an   $n$-bracket,
$n>2$,  on  a  finite  dimensional  vector  space  $V$  such  that
$[y_1,\dots,y_{n-1}]_x=[x,y_1,\dots,y_{n-1}]$   is   a    Filippov
$(n-1)$-bracket for every $x\in V$. Does it imply that
$[\cdot,\dots,\cdot]$ is Filippov itself?}

\medskip\noindent
{\bf Acknowledgment.} The authors are  grateful   to   Jean   Paul
Dufour  and Alexandre Vinogradov for their  remarks  and  comments
and for the referee,  who  pointed  out  the  correct  version  of
(\ref{xx}) and problems concerning the case $n=3$.

\end{document}